\documentclass[11pt, letterpaper]{amsart}

\usepackage{amsmath, amssymb, amsthm} 
\usepackage{graphicx} 
\usepackage[english]{babel}
\usepackage{csquotes}
\tolerance=2000       \usepackage{microtype}
\usepackage[margin=1in]{geometry} 
\usepackage[
  colorlinks=true,
  linkcolor=blue,
  citecolor=blue,
  urlcolor=blue
]{hyperref} 
\usepackage{xcolor} 
\usepackage[
    backend=biber, 
    style=numeric,
    isbn=false,    
    doi=false,    
    url=false,
    maxnames=99
]{biblatex}
\renewbibmacro{in:}{}
\usepackage{orcidlink}

\theoremstyle{plain}
\newtheorem{theorem}{Theorem}[section]
\newtheorem{proposition}[theorem]{Proposition}
\newtheorem{lemma}[theorem]{Lemma}
\newtheorem{corollary}[theorem]{Corollary}
\theoremstyle{definition}
\newtheorem{definition}[theorem]{Definition}
\theoremstyle{remark}
\newtheorem{remark}[theorem]{Remark}

\newcommand{\R}{\mathbb{R}}
\newcommand{\abs}[1]{\left| #1 \right|}
\newcommand{\op}{\operatorname}
\newcommand{\norm}[1]{\left\Vert #1\right\Vert}
\newcommand{\pint}[1]{ \langle #1 \rangle}
\newcommand{\set}[1]{ \lbrace #1 \rbrace}
\newcommand{\intRp}{\int_{\mathbb{R}^+}}

\usepackage{bbm}

\newcommand{\I}[1]{{\mathbbm{1}}_{\set{#1}}}

\usepackage{etoolbox}
\makeatletter
\patchcmd{\@setaddresses}{\indent\emailaddrname}{\emailaddrname}{}{}
\makeatother

\title{Existence and uniqueness for a coagulation-transport equation}

\author[E. Aguilar Anzures]{Edwin Aguilar Anzures \orcidlink{0009-0004-0139-4977}}
\address{
\newline\rm Edwin Aguilar Anzures, L. F. López Ríos and Arno Siri-Jégousse.
\vspace{2mm}
\newline
Instituto de Investigaciones en Matemáticas Aplicadas y en Sistemas \newline
Universidad Nacional Autónoma de  México \newline Circuito Escolar s/n, Ciudad de México C.P. 04510 (México) \vspace{2mm}\newline\textit{Email addresses:} \texttt{luis.aguilar@delta.iimas.unam.mx}, \texttt{luis.lopez@mym.iimas.unam.mx}, \texttt{arno@sigma.iimas.unam.mx}.}

\author[L. F. López Ríos]{Luis F. López Ríos \orcidlink{0000-0001-9949-2620}}

\author[A. Siri-Jégousse]{Arno Siri-Jégousse \orcidlink{0000-0003-1289-2950}}

\begin{document}

\begin{abstract}
    \noindent
    This paper establishes the existence, uniqueness, and regularity of weak solutions to a class of nonlinear coagulation-transport equations describing the asymptotic dynamics of multispecies coalescent processes. The model captures the multi-scale ancestral lineage decay governed by an autonomous drift field, coupled with a time-inhomogeneous quadratic coagulation term parameterized by the initial population mass. By utilizing the metric structure of probability measures, we demonstrate that any weak solution satisfies a unique mild representation formula dictated by the transport semigroup. This characterization allows us to construct a unified narrative where existence, uniqueness, and moment preservation are proved directly through a contractive fixed-point argument in the space of finite Radon measures.
\end{abstract}
\maketitle
\section{Introduction}
Let $\mathcal{M}(\R^+)$ denote the space of finite signed Radon measures and
let $\mathcal{M}_P(\R^+)$ denote its subset of probability measures on $\R^+$. 
Consider the following equation in $\mathcal{M}(\R^+)$,
\begin{equation}\label{eq:coagulationTransport}
    \begin{cases}\displaystyle
\partial_t\mu_t=\partial_x(\psi \mu_t)+\frac{1}{t+\kappa}(\mu_t\star\mu_t-\mu_t),\\
        \mu_0=\nu,
    \end{cases}
\end{equation}
where $\psi:\R^+\to\R^+$ is continuous, $\kappa >0$ and $\nu\in \mathcal{M}_P(\R^+)$ are given parameters. 

In population genetics, this equation arises in the setting of structured genealogical trees.
It describes the asymptotic dynamics of a random multispecies coalescent process with a large but finite number of initial species.
More precisely, it represents the masses of different species genealogies when
the amount of lineages decreases at rate  $\dot x =- \psi(x)$ within each species and $\kappa>0$ is the inverse of the initial population mass, see \cite{lambert2018coagulationtransport, airamTreesWithinTrees}.

This equation may also be obtained by rescaling the Smoluchowski equation \parencite{smoluchowski1917versuch} with mass loss through the division  of the total mass of the system. The Smoluchowski equation with mass loss has been studied in a partial differential equations framework, see \textcite{BanasiakCoagulationMassLoss} and the references therein.

 For a finite signed measure $\mu$ and  $f\in C(\R^+)$ we use the notation $ \pint{\mu,f}=\intRp f(x)\,\mu(dx)$. The space $\mathcal{M}(\R^+)$ is a Banach space with the total variation norm, which we will denote by $\norm{\cdot}$. This space is isometrically isomorphic to the dual space of the continuous functions in $\R^+$ vanishing at infinity, the latter denoted by  $C_0(\R^+)$. 
 Let $C_b(\R^+)$ be the set of continuous bounded functions, $(\mathcal{M}(\R^+),w)$ denotes the set of signed Radon measures on $\R^+$ endowed with the weak topology, i.e., the smallest topology making the map $\mu\mapsto\pint{\mu,f}$ continuous for every $f\in C_b(\R^+)$. Note that, given $T\in(0,\infty]$, a map $\mu_t:[0,T]\to(\mathcal{M}(\R^+),w)$ is continuous if and only if the function $t\mapsto\pint{\mu_t,f}$ is continuous in $[0,T]$ for every $f\in C_b(\R^+)$. We will refer to this property equally as weak continuity. We will see that the weak topology is the natural topology for studying the transport operator. 

 Since we are interested in a weak (measure-valued) formulation of this equation, we will take into consideration for our study the natural test space
\[
\mathcal{E}=\set{f\in C^1(\R^+)\,:\,f, \psi f' \in C_b(\R^+)}.
\]

\begin{definition}[Weak solution] Let $\kappa > 0$ and $\nu\in \mathcal{M}_P(\R^+)$. A probability-valued process $(\mu_t)_{t\geq0}$ is a weak solution of  \eqref{eq:coagulationTransport} with initial condition $\nu $  if the following conditions are satisfied:
\begin{enumerate}
    \item[$(i)$] The map $t\mapsto\mu_t$ is weakly continuous;
\item[$(ii)$] For every $t\geq 0$,
\begin{equation}
    \label{eq:weakSolution}
    \pint{\mu_t,f}=\pint{\nu,f}-\int_0^t\pint{\mu_s,\psi f'}\,ds+\int_0^t\frac{1}{s+\kappa}\left[\pint{\mu_s\star \mu_s, f}-\pint{\mu_s, f}\right]\,ds\qquad \text{for all }f\in\mathcal{E} .
\end{equation}
\end{enumerate}
\end{definition} 

The problem of finding weak solutions to \eqref{eq:coagulationTransport} first appeared in \textcite{lambert2018coagulationtransport} where they obtained an existence and uniqueness result by means of probabilistic methods based on branching processes. Similar weak formulations have appeared in \cite{Norris1999} and \cite{Fournier2006}.

More precisely, \textcite{Norris1999} established the general framework to treat weak formulations for Smoluchowski type equations, that is, evolution of measures performing only pure coagulation phenomena. There, the author uses measure-theoretic techniques for non-negative mass distributions. The present work follows a more abstract functional-analytic perspective. We utilize the Riesz-Markov-Kakutani Representation Theorem \cite[Theorem 14.14]{folland1999real} to define the evolution of the system via Pettis integration \parencite{pettis1938integration}. This approach allows for the inclusion of  finite signed Radon measures and provides a robust way to prove existence and uniqueness of weak solutions for  equations describing transport and coagulation and, as a consequence, extending to a wider range of functions $\psi$ the already obtained results. 
This is our main result, which we will prove at the end of Section \ref{sec:MildFormulation}.

\begin{theorem}\label{thm:mainTheorem}
Let $\psi$ be a non-negative locally Lipschitz function and let $\nu$ be a probability measure on $\R^+$. There exists a unique weak solution of \eqref{eq:coagulationTransport} with initial condition $\nu$.  
\end{theorem}

In the case that the initial condition meets certain regularity conditions, for instance being absolutely continuous with respect to Lebesgue measure on $\R^+$ or it is a Dirac measure, we can prove that the solution preserves, at least, these properties; for all positive time the solution remains absolutely continuous or keeps always a singular part plus an absolutely continuous part respectively. In this regard, we prove the following theorem in Section \ref{sec:Regularity},
\begin{theorem} \label{thm:secondMainThm}
Let $\mu_t$ be a weak solution to \eqref{eq:coagulationTransport} with initial condition $\nu\in \mathcal{M}_P(\R^+)$. Then,
\begin{itemize}
\item[$(i)$] If $\nu$ has finite first moment, then $\mu_t$ has finite first moment for all $t>0$.
    \item[$(ii)$] If $\nu\ll\mathcal{L}$, then $\mu_t\ll\mathcal{L}$ for all $t>0$.
    \item[$(iii)$] If $\nu=\boldsymbol{\delta}_{x_0}$ for some $x_0\in\R^+$ and $\psi(x)=x^2/2$, then there exists a function $u:[0,T]\times \R^+\to\R^+$ such that $u(t,\cdot)\in L^1(\R^+)$ for all $t\in[0,T]$ and 
\begin{equation*}
\mu_t(dx)=\frac{\kappa}{t+\kappa }\boldsymbol{\delta}_{\frac{2x_0}{2+tx_0}}+u(t,x)dx.
\end{equation*}
\end{itemize}
\end{theorem}
This last case describes the the asymptotic dynamics of a nested Kingman's process \cite{airamNestedCoalescent} with a large but finite initial number of species, all with the same number of lineages each. In this case each species starts with a proportion of $x_0$ lineages-mass, while the mass within each species evolves as $\dot x=-x^2/2$. The proof of this theorem is essentially provided in Propositions \ref{prop:lebesgueRegularity} and \ref{prop:DeltaSolutionFormula}.

From the biological point of view the natural space to frame this problem is $\mathcal{M}_P(\R^+)$, reflecting the probability of observing species with a given amount of rescaled gene lineages. Even so, $\mathcal{M}_P(\R^+)$ lacks the Banach space structure, preventing a direct application of standard fixed point theorems. To overcome this structural limitation, we extend the boundaries of our search space and formulate the well-posedness problem within the Banach space of finite signed Radon measures. This extended setting allows us to apply properly the Banach fixed point theorem to obtain local existence and uniqueness. Subsequently, we will prove that the fixed points whose initial condition is a probability measure remains in $\mathcal{M}_P(\R^+)$, thereby recovering the biological probability interpretation a posteriori.

\begin{definition}(\textbf{Local solution})\label{def:localSolution}
Let $T > 0$. A weakly continuous map $\mu_t:[0,T]\to \mathcal{M}(\R^+)$ is a (weak) local solution of  \eqref{eq:coagulationTransport} with initial condition $\nu\in \mathcal{M}(\R^+)$  if, \eqref{eq:weakSolution} holds for every $t\in[0,T]$.    
\end{definition}

We can derive the probability mass preservation  ipso facto from the previous definition.

\begin{lemma}[\textbf{Probability conservation}]\label{lem:probabilityConcervation}
    Let $\mu_t:[0,T]\to \mathcal{M}(\R^+)$ be a local weak solution  with initial condition $\nu$. If $\nu(\R^+)=1$, then $\mu_t(\R^+)=1$ for $t\in[0,T]$.
\end{lemma}
\begin{proof}
By taking $f=1$ in  \eqref{eq:weakSolution} we obtain 
\[
\mu_t(\R^+)=1+\int_0^t\frac{1}{s+\kappa}[\mu_s(\R^+)^2-\mu_s(\R^+)]\,ds
\]
Since the map $t\mapsto \mu_t(\R^+)$ is continuous, from the integral equation we deduce that it is also absolutely continuous, and we can differentiate almost everywhere. We obtain the following identity: 
\[
\frac{d}{dt}\mu_t(\R^+)=\frac{1}{t+\kappa}[\mu_t(\R^+)^2-\mu_t(\R^+)].
\]
It turns out that the only solution with $\mu_0(\R^+)=1$ is $\mu_t(\R^+)=1$ for all $t>0$. 
\end{proof}

\section{Mild formulation and fixed point map}\label{sec:MildFormulation}
The goal of this section is to introduce an appropriate space where we can rely on the Banach fixed point theorem to find a mild solution on a possibly small time interval. Also we will show that the mild solution or, more properly, the fixed point of a suitable operator, 
is also a local weak solution. 
Later we will see that we can extend the obtained solution to an arbitrary time by using a standard argument and the fact that the local weak solutions preserve the positivity and the mass.

We will develop the analysis in the space $\mathcal{M}(\R^+)$ selecting an appropriate topology. 
Since we pretend to apply a fixed point argument, we need an underlying Banach space structure with a norm compatible with the weak continuity of paths and strong enough to handle the nonlinearity of the coagulation. 

We will denote a map $\mu_t:[0,T]\to\mathcal{M}(\R^+)$ simply by $\mu_t$. We will also consider the space $C([0,T];(\mathcal{M}(\R^+),w))$ of the weakly continuous maps
and the space $L^\infty([0,T];\mathcal{M}(\R^+))$ of bounded maps, this is a Banach space with the norm  
\[
\norm{\mu_t}_{\infty}:=\sup_{0< t< T}\norm{\mu_t}.
\] 

\subsection{Transport of measures}\label{subSec:tranportOfMeasures}

We borrow the notation and terminology of the optimal transport literature; see for  instance \textcite{AmbrosioGigliSavare2008}.
\begin{definition}[\textbf{Push-forward of a measure}]
    If $\mu$ is a finite signed measure on $\R^+$ and $\Phi:\R^+\to\R^+$ is a Borel map, we denote by $\Phi_\#\mu$ the push-forward of $\mu$ through $\Phi$, defined by 
    \[
\Phi_\#\mu(B):=\mu(\Phi^{-1}(B)),\quad\text{for all } B\in\mathcal{B}(\R^+).
    \]
\end{definition}

We are interested in whether  $\psi$ can generate a unique flow $\Phi_t$, that is, the one solving the equation 
\begin{equation}\label{eq:flow}
    \frac{d}{dt}\Phi_t(x) = -\psi(\Phi_t(x)),
\quad \Phi_0(x)=x
\end{equation}
for each fixed $x\in\R^+$. If the flow exists we will see that we can construct a weak solution to the coagulation transport equation by  using the push-forward measure of the flow as our transport operator. Notwithstanding, the uniqueness of the solution strongly depends on the uniqueness of the flow. Since the uniqueness of the flow is of such relevance, we can play with more general functions $\psi$, even with less regularity, and the arguments will still be valid. For an extensive treatment on whether the ODE can have a solution with low regularity parameters  we refer the reader to \textcite{DiPerna1989} and \textcite{CrippaDeLellis2008}.

Let us assume that $\psi:\R^+\to\R^+$ is continuous and such that the equation 
\eqref{eq:flow} has a unique $C^1$ solution that satisfies the semigroup properties
\begin{equation}
\label{eq:semigroupProperties}  \Phi_{s+t}=\Phi_{s}\circ\Phi_{t}\quad\text{and}\quad\Phi_0=\op{Id}.
\end{equation}
 For instance, it is possible to achieve this by assuming $\psi$ to be locally Lipschitz. Let us  denote $
S_{\Phi}(t)\mu := (\Phi_{t})_\#\mu,$ for $\mu\in\mathcal{M}(\R^+)$. We will call it the transport operator.

The transport operator preserves the total variation norm, i.e., for every $t\geq 0$ 
\begin{equation}
    \label{eq:transporBound}
\norm{S_\Phi(t)\mu}=\norm{\mu},\quad\text{for all}\quad\mu\in\mathcal{M}(\R^+).
\end{equation} The transport operator just moves mass, it does not create it nor destroy it. Thus, this operator leaves the cone of positive measures invariant in $(\mathcal{M}(\R^+),\norm{\cdot})$. However, the semigroup $(S_{\Phi}(t))_{t\geq0}$ is not absolutely continuous with respect to the total variation norm $\norm{\cdot}$ since $S_{\Phi}(t)\boldsymbol{\delta}_x=\boldsymbol{\delta}_{\Phi_t(x)}$ and $\norm{\boldsymbol{\delta}_x-\boldsymbol{\delta}_y}=2$ whenever $x\neq y$. As a consequence, $(\mathcal{M}(\R^+),\norm{\cdot})$ is not a suitable space on its own to study the transport operator, we will need to impose certain conditions guaranteeing the weak convergence of the measures.

The transport operator satisfies weakly the transport equation, as it were. This result is standard and the real difficulty is justifying the differentiation under the integral sign. 
\begin{proposition}
\label{prop:weakContinuityTransport}
Let $\nu\in\mathcal{M}(\R^+)$. If $\psi$ is locally Lipschitz, then for every test function $f\in \mathcal{E}$, the map
\[
t\mapsto\pint{S_{\Phi}(t)\nu,f}
\]
is $C^1$ and satisfies 
\[
\frac{d}{dt} \pint{S_{\Phi}(t)\nu,f}=-\pint{S_{\Phi}(t)\nu,\psi f'}.
\] That is, $(S_{\Phi}(t)\nu)_{t\geq0}$ is a weak solution of the transport equation
\[
\begin{cases}
    \partial_t \mu_t =\partial_x(\psi\mu_t),\quad t\geq0,\\
    \mu_0=\nu.
\end{cases}
\]
\end{proposition}
\begin{proof}
By definition of the push-forward of a measure, we have 
\[
\pint{S_{\Phi}(t)\nu,f} = \intRp f(\Phi_t(x))\,\nu(dx).
\] The map $t\mapsto f(\Phi_t(x))$ is $C^1$, bounded and its derivative $t\mapsto-f'\psi(\Phi_t(x))$ is also bounded, we  can differentiate under the integral sign so that,
\[
\frac{d}{dt} \intRp f(\Phi_t(x))\,\nu(dx)=-\intRp f'(\Phi_t(x))\psi (\Phi_t(x))\,\nu(dx).
\]
Thus, 
\[
\frac{d}{dt}\pint{S_\Phi(t)\nu,f}=-\pint{S_\Phi(t)\nu,f'\psi}.
\]
\end{proof}
\begin{remark}
As a consequence of this result we have that, when $\nu\in\mathcal{M}(\R^+)$, the map $t\mapsto S_\Phi(t)\nu$ is an element of  $C([0,T];(\mathcal{M}(\R^+),w))\cap L^\infty([0,T];\mathcal{M}(\R^+))$.
\end{remark}

\subsection{Coagulation operator}
 Let us consider the coagulation operator $Q:\R^+\times \mathcal{M}(\R^+)\to \mathcal{M}(\R^+)$ defined by 
\[
Q(t,\mu)=\frac{1}{t+\kappa}(\mu\star\mu-\mu).
\]
This operator  is well defined since $(\mathcal{M}(\R^+),\norm{\cdot})$ is a Banach algebra with the convolution, in the sense that the convolution operator acts as a product in this space: $\mathcal{M}(\R^+)$ is closed under convolutions and the operation is sub-multiplicative, that is, $\norm{\nu\star\mu}\leq \norm{\nu}\norm{\mu}$ for every $\mu,\nu\in\mathcal{M}(\R^+)$. The last inequality is known as Young's convolution inequality.

\begin{lemma}\label{lem:cobo}
    For all $\mu, \tilde{\mu}\in \mathcal{M}(\R^+) $, the following estimation holds
    \[
    \norm{\mu \star \mu-\tilde{\mu} \star \tilde{\mu}}\leq (\norm{\mu}+\norm{\tilde{\mu}})\norm{\mu - \tilde{\mu}}.
    \]
\end{lemma}
\begin{proof}
   We only need to apply that the convolution is distributive over addition:
    \begin{align*}
        \norm{\mu \star \mu-\tilde{\mu} \star \tilde{\mu}} &=\norm{\mu \star \mu-\tilde{\mu}\star \mu+\tilde{\mu}\star \mu-\tilde{\mu} \star \tilde{\mu}}\\
        &=\norm{(\mu - \tilde{\mu}) \star \mu + \tilde{\mu} \star (\mu - \tilde{\mu})}\\
        &\leq \norm{(\mu - \tilde{\mu}) \star \mu}+\norm{\tilde{\mu} \star (\mu - \tilde{\mu})}\\
        &\leq (\norm{\mu}+\norm{\tilde{\mu}})\norm{\mu - \tilde{\mu}}.
    \end{align*}
\end{proof}

This seemingly trivial inequality has many revealing consequences. For instance, the following corollary shows that the operator $Q$ is quadratic and locally Lipschitz in an appropriate space of maps $t\mapsto\mu_t$.
\begin{corollary}[Coagulation bound]
    Given $\mu, \tilde{\mu}\in \mathcal{M}(\R^+) $, the following bound holds for every $t\geq0$:
    \begin{equation}\label{eq:CoagulationLipschitzBound}
        \norm{Q(t,\mu)-Q(t,\tilde{\mu})}\leq \frac{1}{t+\kappa}\norm{\mu-\tilde{\mu}}(\norm{\mu}+\norm{\tilde{\mu}}+1).
    \end{equation}
\end{corollary}
\begin{proof}

Consider two measures $\mu$ and $\tilde{\mu}$. Then
\begin{align*}
    \norm{Q(t, \mu) - Q(t, \tilde{\mu})} &= \frac{1}{t + \kappa}\norm{\mu \star \mu-\mu-\tilde{\mu} \star \tilde{\mu}+ \tilde{\mu}}\\
    &\leq \frac{1}{t + \kappa}(\norm{\mu \star \mu-\tilde{\mu} \star \tilde{\mu}}+\norm{\mu-\tilde{\mu}}).
\end{align*}
Applying Lemma \ref{lem:cobo}, we obtain the desired bound.
\end{proof}

Now, let us focus on the regularity implications.
\begin{corollary}
The operator $Q$ has the following properties:
\begin{itemize}
    \item[$(i)$]  If $\mu_t\in L^\infty([0,T];\mathcal{M}(\R^+))$, then $Q(t,\mu_t)\in L^\infty([0,T];\mathcal{M}(\R^+))$. 

    \item[$(ii)$] The map $Q$ from $L^\infty([0,T];\mathcal{M}(\R^+))$ to itself defined as  $\mu_t\mapsto Q(t,\mu_t)$ is locally Lipschitz and thus continuous.
\end{itemize}
\end{corollary}
\begin{proof} Applying \eqref{eq:CoagulationLipschitzBound} we obtain
\begin{equation} \label{eq:coagulationBound}
    \norm{Q(t,\mu)}\leq \frac{\norm{\mu}(\norm{\mu}+1)}{t+\kappa}\quad\text{for every }t\in [0,T].
\end{equation} Thus, 
\[
\norm{Q(t,\mu_t)}\leq\frac{1}{\kappa}\norm{\mu_t}_\infty(\norm{\mu_t}_\infty+1)<\infty\quad\text{for every }t\in[0,T],
\] thereby $\norm{(Q(t,\mu_t))}_\infty<\infty$. This shows $(i)$.

Now, let us consider $\mu_t,\tilde\mu_t\in L^\infty([0,T];\mathcal{M}(\R^+))$. 
Applying \eqref{eq:CoagulationLipschitzBound} again and taking the supremum over $t\in[0,T]$ we obtain 
\[
\norm{Q(t,\mu_t)-Q(t,\tilde\mu_t)}_\infty\leq \frac{1}{\kappa}\norm{\mu_t-\tilde{\mu}_t}_\infty(\norm{\mu_t}_\infty+\norm{\tilde{\mu}_t}_\infty+1).
\]
Thus, if $\tilde\mu_t$ happens to be in a ball centered at $\mu_t$, this bound gives us the local Lipschitz condition. This proves $(ii)$.
\end{proof}

\begin{proposition}
\label{prop:weakContinuityConvolution}
    If $\mu_t\in C([0,T];(\mathcal{M}(\R^+),w))\cap L^\infty([0,T];\mathcal{M}(\R^+))$, then $\mu_t\star\mu_t$ is weakly continuous.
\end{proposition}

\begin{proof}
 Let us fix $\tilde t\in[0,T]$, a sequence $(t_n)_{n\geq 0}$ in $[0,T]$ such that $t_n\to \tilde t$ as $n\to\infty$ and a test function $f\in C_b(\R^+)$ during the proof. For $t\in[0,T]$ and $x\in\R^+$ let us define 
 \[
 \varphi_t(x):=\intRp f(x+y)\, \mu_t(dy).
 \] 
 We will use this adaptation to prove the weak convergence result. Note that we have the bound $ \abs{\varphi_{t}(x)}\leq \norm{f}_\infty\norm{\mu_t}_\infty$, that is, $(t,x)\mapsto\varphi_t(x)$  is uniformly bounded. Moreover, since $f$ is continuous and bounded and $\mu_t$ is bounded in total variation, dominated convergence theorem implies that for each $t\in[0,T]$, the map $x\mapsto \varphi_t(x)$ is bounded and continuous, so $\pint{\mu_t\star\mu_t,f}=\pint{\mu_t,\varphi_t}$.

Let $n\geq0$, by adding and subtracting $\pint{\mu_{t_n},\varphi_{\tilde t}}$ we get
 \begin{equation}
\label{eq:convolutionWeakBound}
\abs{\pint{\mu_{t_n}\star\mu_{t_n},f}-\pint{\mu_{\tilde t}\star\mu_{\tilde t},f}}\leq \abs{\pint{\mu_{t_n},\varphi_{t_n}-\varphi_{\tilde t}}}+\abs{\pint{\mu_{t_n},\varphi_{\tilde t}}-\pint{\mu_{\tilde t},\varphi_{\tilde t}}}.
 \end{equation}
Since $\varphi_{\tilde t}\in C_b(\R^+)$, using weak convergence of $\mu_{t_n}$ to $\mu_{\tilde t}$ we have $\abs{\pint{\mu_{t_n},\varphi_{\tilde t}}-\pint{\mu_{\tilde t},\varphi_{\tilde t}}}\xrightarrow[n\to\infty]{}0$. It remains to argue why $\abs{\pint{\mu_{t_n},\varphi_{t_n}-\varphi_{\tilde t}}}$ goes to $0$ as $n\to\infty$. For this, we will show first that the sequence $(\varphi_{t_n})_{n\geq0}$ converges uniformly to $\varphi_{\tilde t}$ on a compact set and later use a tightness argument to bound the tail when integrating against $\mu_{t_{n}}$.

 Let us fix a compact interval $J\subset\R^+$, we will prove uniform convergence of $(\varphi_{t_n})_{n\geq0}$ over $J$. Since $\R^+$ is a Polish space and $\mu_{t_n}$ converges weakly to $\mu_{\tilde t}$, in this way we can  apply  Prohorov's theorem (\textcite[Theorem 8.6.2]{bogachev2007measure}) and obtaining that the family $\set{\mu_{t_n},\,n\geq0}\cup\set{\mu_{\tilde t}}$ is uniformly tight. Therefore,  for $\eta>0$ we can choose a compact interval $K\subset\R^+$ such that for all $n\geq0$, $\abs{\mu_{t_n}}(\R^+\setminus K)<\eta$; then, for  $x,\tilde x\in J$,
\begin{align*}
\abs{\varphi_{_{t_n}}(x)-\varphi_{_{t_n}}(\tilde x)}&\leq \int_{K} \abs{f(x+y)-f(\tilde x-y)}\abs{\mu_{t_n}}(dy)+\int_{\R^+\setminus K}\abs{f(x+y)-f(\tilde x-y)}\abs{\mu_{t_n}}(dy)\\
    &\leq \int_{K} \abs{f(x+y)-f(\tilde x-y)}\abs{\mu_{t_n}}(dy) + 2\norm{f}_\infty \norm{\mu_{t_n}}_\infty\eta.
\end{align*}
To finish the argument, let us note that $f$ is uniformly continuous in the compact interval $J$, so we can choose $\kappa>0$ such that for any $z,\tilde z \in J$, $\abs{f(z)-f(\tilde z)}<\eta$ whenever $\abs{z-\tilde z}<\kappa$. 
As a consequence, if $\abs{x-\tilde x}<\kappa$ we can bound the remaining integral as
\begin{align*}
    \int_{K} \abs{f(x+y)-f(\tilde x-y)}\abs{\mu_{t_n}}(dy)\leq \eta \abs{
    \mu_{t_n}}(K)\\
    \leq \eta \norm{
    \mu_{t_n}}_\infty.
\end{align*}
 This implies  $\abs{\varphi_{{t_n}}(x)-\varphi_{{t_n}}(\tilde x)}\leq \eta \norm{\mu_t}_\infty(1+2\norm{f}_\infty)$, meaning that   the sequence $(\varphi_{t_n})_{n\geq 0}$ is equicontinuous on $J$.  
Moreover, the Arzelà-Ascoli theorem guarantees that the sequence has a uniformly convergent subsequence. By the weak continuity, we have that $\varphi_{t_n}(x)\xrightarrow[n\to\infty]{}\varphi_{\tilde t}(x)$, that is, $\varphi_{t_n}$ converges point-wise to $\varphi_t$, thus the limit of such subsequence must be exactly $\varphi_{\tilde t}$. 
Since we can apply the same argument to every subsequence of $\set{\varphi_{t_n}}_{n\geq 0}$ we conclude that actually $(\varphi_{t_n})_{n\geq 0}$ converges uniformly to   $\varphi_{\tilde t}$ over the compact inverval $J$.

Finally, let $\varepsilon>0$. We can choose again a compact interval $I$ such that for all $t\in\set{t_n\ ,\ n\geq0 }\cup\set{\tilde t}$, $\abs{\mu_t}(\R^+\setminus I)<\varepsilon$. Also, using that $(\varphi_{t_n})_{n\geq0}$ converges uniformly to $\varphi_{\tilde t}$ on the compact interval $I$ we can choose $n\in\mathbb{N}$ big enough such that $\sup_{x\in I} \abs{\varphi_{t_n}(x)-\varphi_{\tilde t}(x)}<\varepsilon$. Then,
\begin{align*}
\abs{\pint{\mu_{t_n},\varphi_{t_n}-\varphi_{\tilde t}}} &\leq \int_I\abs{\varphi_{t_n}(x)-\varphi_{\tilde t}(x)}\, \abs{\mu_{t_n}}(dx)+\int_{\R^+\setminus I}\abs{\varphi_{t_n}(x)-\varphi_{\tilde t}(x)}\, \abs{\mu_{t_n}}(dx)\\
&\leq \varepsilon \norm{\mu_t}_\infty+2\norm{f}_\infty\norm{(\mu_t)}_\infty\varepsilon.
\end{align*}
This shows that $\abs{\pint{\mu_{t_n},\varphi_{t_n}-\varphi_{\tilde t}}}\xrightarrow[n\to\infty]{}0$. Then, taking the limit in \eqref{eq:convolutionWeakBound} we obtain the desired result.
\end{proof}

\begin{remark}
    This previous result implies that, whenever $\mu_t\in C([0,T];(\mathcal{M}(\R^+),w))$, the map $t\mapsto Q(t,\mu_t)$ will be weakly continuous as well.
\end{remark}

\subsection{Mild formulation} In this section we construct an evolution operator in a suitable subspace of $C([0,\infty);(\mathcal{M}(\R^+),w))$ which will ultimately allow us to prove the existence of solutions via a fixed-point argument. To achieve this we extend the Duhamel's formula to the weak topology of measures. More specifically,  the operator will be defined as the sum of the free evolution of the initial condition under the transport operator plus the weak time-integral of the transport operator applied to the coagulation operator.

First focus on finding an appropriate space with Banach structure to apply the fixed point theorem: As we saw in the previous sections, weak continuity is the natural concept for transport problems, since total variation norm is too strong for the continuity of the transport. Notwithstanding, the total variation norm provides the correct estimates for contractions. So, let us consider the space 
\[
X_T=C([0,T];(\mathcal{M}(\R^+),w))\cap L^\infty([0,T],\mathcal{M}(\R^+))
\]
\begin{proposition}
    The space $(X_T,\norm{\cdot})$ is Banach.
\end{proposition}
\begin{proof}
Let $(\mu_t^n)$ be a Cauchy sequence in $X_T$. 
There exists $\mu_t\in L^\infty([0,T],\mathcal{M}(\R^+))$ such that 
\[
\sup_{0< t< T}\norm{\mu^n_t-\mu_t}\xrightarrow[n\to\infty]{}0.
\]
Let us fix $f\in C_b(\R^+)$.
Since $\abs{\pint{\mu^n_t,f}-\pint{\mu_t,f}}\leq\norm{f}_{L^\infty(\R^+)}\norm{\mu_t^n-\mu_t}$, we have 
\[
\sup_{0\leq t\leq T} \abs{\pint{\mu_t^n,f}-\pint{\mu_t,f}}\xrightarrow[n\to\infty]{}0.
\]
This means that the functions $t\mapsto\pint{\mu^n_t,f}$ converge uniformly on $[0,T]$ to $t\mapsto\pint{\mu_t,f}$.
Since that uniform limit of continuous functions is also continuous, we conclude that the map is continuous. This implies that $\mu_t\in X_T$. 
\end{proof}

Duhamel's formula involves integrating on a trajectory on a Banach space. In this case to rigorously define the Pettis integral for a measure-valued trajectory a crucial step is verifying that the scalar function $t\mapsto\norm{\mu_t}$ is Borel measurable on $[0,T]$ whenever $\mu_t\in X_T$. We can ensure this by showing that the map is lower semicontinuous.
    To see this we will use the duality relationship $\mathcal{M}(\R^+)\cong C_0(\R^+)^*$ given by the Riesz-Markov-Kakutani Representation Theorem (\textcite[Theorem 7.17]{folland1999real}). Then, for $\mu\in \mathcal{M}(\R^+)$ we can compute its $\mathcal{M}(\R^+)$ norm as follows; 
\begin{equation}
\label{eq:normDuality}\norm{\mu}=\sup\set{\abs{\pint{\mu,f}}\, :\, f\in C_0(\R^+)}.
\end{equation}
Let us then pick a countable dense subset $\set{f_n}_n$ of $C_0(\R^+)$, it is well known that this space is separable with the supremum norm for continuous functions. 
Then, can compute the norm of a measure $\mu\in\mathcal{M}(\R^+)$ as $\norm{\mu}=\sup_{n\in\mathbb{N}}\abs{\pint{\mu,f_n}}$. 
Now, since the map $t\mapsto \pint{\mu,f_n}$ is continuous and $t\mapsto\norm{\mu_t}$ is the supremum of continuous functions, it follows that $t\mapsto\norm{\mu_t}$ is lower semicontinuous and indeed Borel measurable.

We are now equipped with the necessary framework to construct an evolution  operator in $X_T$ to apply a fixed point argument. Let $\mu_t\in X_T$ and $f\in C_b(\R^+)$. Then, Propositions \ref{prop:weakContinuityTransport} and \ref{prop:weakContinuityConvolution} assure that for a fixed $t\in[0,T]$ the map $s\mapsto \pint{S_{\Phi}(t-s)\,Q(s,\mu_s),f}$ is continuous for $s\in[0,T]$, so it can be integrated from $0$ to $t$; also applying the bounds \eqref{eq:coagulationBound} and \eqref{eq:transporBound} we see that $\abs{\pint{S_{\Phi}(t-s)\,Q(s,\mu_s),f}}\leq \frac{1}{\kappa}\norm{f}_\infty\norm{\mu_t}_\infty(1+\norm{\mu_t}_\infty)$. The duality $\mathcal{M}(\R^+)\cong C_0(\R^+)^*$ allows us to define $\mathcal{I}(\mu_t)\in \mathcal{M}(\R^+)$ trough its action in $C_0(\R^+)$ as 
    \begin{equation}
    \label{eq:pettisIntegralDefinition}
        \pint{\mathcal{I}(\mu_t),f}:=\int_0^t\pint{S_{\Phi}(t-s)\,Q(s,\mu_s),f}\,ds\quad\text{for all}\quad f\in C_0(\R^+).
    \end{equation}
 This procedure to define $\mathcal{I}$ using duality is known as the Pettis integral in $\mathcal{M}(\R^+).$ 

\begin{proposition}
\label{prop:pettisIntegralBound}
     If $\mu_t\in X_T$, then
     $\mathcal{I}(\mu_t)\in X_T$ and 
\begin{equation}\label{eq:pettisIntegralBound}
         \norm{\mathcal{I}(\mu_t)}_\infty\leq (\norm{\mu_t}_\infty+\norm{\mu_t}_\infty^2)\log\left(\frac{T+\kappa}{\kappa}\right).
     \end{equation}
\end{proposition}
\begin{proof}
To obtain the bound we will compute the norm by taking advantage of the duality relationship \eqref{eq:normDuality}. Let $f\in C_0(\R^+)$. Taking absolute value on the definition \eqref{eq:pettisIntegralDefinition} and using that $S_{\Phi}(t-s)\,Q(s,\mu_s)\in \mathcal{M}(\R^+)$, for every $0\leq s<t< T$, we get, for every $t\in[0,T]$,
\begin{align*}
    \abs{\pint{\mathcal{I}(\mu_t),f}}&\leq \int_0^t\intRp\abs{f(x)}\, \abs{S_{\Phi}(t-s)\,Q(s,\mu_s)}(dx)\ ds\\
    &\leq \int_0^t\intRp(1+x)\, \abs{S_{\Phi}(t-s)\,Q(s,\mu_s)}(dx)\ ds
    =\int_0^t\norm{S_{\Phi}(t-s)\,Q(s,\mu_s)} \ ds.
\end{align*}
Thus, applying \eqref{eq:coagulationBound} and \eqref{eq:transporBound} we get \[
    \abs{\pint{\mathcal{I}(\mu_t),f}}\leq (\norm{\mu_t}_\infty+\norm{\mu_t}_\infty^2)\int_0^T
\frac{1}{s+\kappa}\,ds. 
\] 
Taking the supremum over all  $f\in C_0(\R^+)$, we obtain the desired inequality \eqref{eq:pettisIntegralBound}, which also proves that $\mathcal{I}(\mu_t)\in L^\infty([0,T],\mathcal{M}(\R^+))$. 

It remains to show that $\mathcal{I}(\mu_t)$ is weakly continuous. Let $f\in C_0(\R^+)$ and fix $t\in[0,T]$. For $h>0$, 
\begin{align*}
\abs{\pint{\mathcal{I}(\mu_{t+h}),f}-\pint{\mathcal{I}(\mu_{t}),f}}\leq& \int_t^{t+h}\abs{\pint{S_{\Phi}(t+h-s)\,Q(s,\mu_s),f}}\, ds\\
&+\int_0^{t}\abs{\pint{S_{\Phi}(t+h-s)\,Q(s,\mu_s)-S_{\Phi}(t-s)\,Q(s,\mu_s),f}}\, ds
\end{align*}
Now observe that we can bound the integral of the first term:
\[
\int_t^{t+h}\abs{\pint{S_{\Phi}(t+h-s)\,Q(s,\mu_s),f}}\, ds\leq \norm{f}_\infty (\norm{\mu_t}_\infty+\norm{\mu_t}_\infty^2)\int_t^{t+h}\frac{1}{s+\kappa}\,ds. 
\]
So the first term goes to $0$ as $h\to0$. 
On the other hand, since the map $h\mapsto S_{\Phi}(t+h-s)\,Q(s,\mu_s)$ is weakly continuous for each $s\in[0,T]$, we know that the integrand converges to zero as $h\to0$.
In other words, the integrand converges point-wise to zero as $h\to0$. And, since the integrand is bounded in the following way:
\[
\abs{\pint{S_{\Phi}(t+h-s)\,Q(s,\mu_s)-S_{\Phi}(t-s)\,Q(s,\mu_s),f}}\leq 2\norm{f}_\infty \norm{\mu_t}_\infty^2)\frac{1}{s+\kappa},
\]
the dominated convergence theorem assures that the second therm converges to zero as $h\to0$. This proves that \[
\abs{\pint{\mathcal{I}(\mu_{t+h}),f}-\pint{\mathcal{I}(\mu_{t}),f}}\xrightarrow[h\to0]{}0,
\]
and, thus, the map $t\mapsto\pint{\mathcal{I}(\mu_{t}),f}$ is continuous for each $f\in C_0(\R^+)$. But, since $t\mapsto\norm{\mathcal{I}(\mu_{t})}$ is uniformly bounded we conclude that indeed $\mathcal{I}(\mu_{t})\in X_T$
This finishes the proof.
\end{proof}

Now we can define the fixed point operator $\mathcal{F}$,  
\begin{equation}\label{eq:mildFormulation}
    \mathcal{F}(\mu_t) := S_{\Phi}(t)\nu+  \mathcal{I}(\mu_t),\quad\text{ for }\mu_t\in X_T,
\end{equation}
where $S_{\Phi}$  is the transport operator defined in Section \ref{subSec:tranportOfMeasures}. Propositions \ref{prop:pettisIntegralBound} and \ref{prop:weakContinuityTransport} guarantee that this is a well-defined map $\mathcal{F}:X_T\to X_T$. Later we will prove that $\mathcal{F}$ maps a closed subset of $X_T$ into itself. We need to verify that a local weak solution is a fixed point of $\mathcal{F}$, and also, that a fixed point of $\mathcal{F}$ satisfies
the weak formulation \eqref{eq:weakSolution}. The next proposition shows that the mild formulation is compatible with the weak formulation. 

\begin{proposition}[Weak implies mild] \label{prop:weakIsMild} If $\mu_t$ is a non-negative local weak solution, then $\mathcal{F}(\mu_t)=\mu_t$ for all $t\in[0,T]$.
\end{proposition}
\begin{proof}
 Let $f\in C^1_c(\R^+)$ and $\tilde t\in[0,T]$ fixed. Since $t\mapsto\pint{\mu_t,f}$ is differentiable also $f\circ\Phi_{\tilde t-t}$ for all $t\in(0,\tilde t\,)$, we can differentiate using the product rule as follows 
\[
\frac{d}{dt}\pint{\mu_t,f\circ\Phi_{\tilde t-t}}=-\pint{\mu_t,\psi(f\circ\Phi_{\tilde t-t})'}+\frac{1}{t+\kappa}(\pint{\mu_t\star\mu_t,f\circ\Phi_{\tilde t-t}}+\pint{\mu_t,f\circ\Phi_{\tilde t-t}})+\pint{\mu_t,\frac{d}{dt}f\circ\Phi_{\tilde t-t}}.
\]
Since the flow $\Phi$ satisfies the identity $\psi \frac{d}{dx}\Phi_t=\psi\circ\Phi_t$ almost everywhere, we see that the transport terms cancel out, leaving 
\[
\frac{d}{dt}\pint{\mu_t,f\circ\Phi_{\tilde t-t}}=\frac{1}{t+\kappa}(\pint{\mu_t\star\mu_t,f\circ\Phi_{\tilde t-t}}+\pint{\mu_t,f\circ\Phi_{\tilde t-t}}),
\]
where we can identify in the right-hand side the pairing  $\pint{Q(t,\mu_t),f\circ\Phi_{\tilde t-t}}$. Integrating from $0$ to $\tilde t$ we get 
\[
\pint{\mu_t,f}-\pint{\nu,f\circ\Phi_t}=\int_0^{\tilde t}\pint{Q(t,\mu_t),f\circ\Phi_{\tilde t-t}}\, dt.
\]
Finally, rearranging the terms and using the definition of the transport operator on measures we get 
\[
\pint{\mu_{\tilde t},f}=\pint{\mathcal{F}\mu_{\tilde t},f}.
\]
Since this is valid for every $f\in C^1_c(\R^+)$ and $\tilde t\in[0,T]$ we conclude that indeed $\mu_t=\mathcal{F}(\mu_{ t})$ for every $t\in[0,T]$, that is, $\mu_t$ is a fixed point of $\mathcal{F}$ in $X_T$. 
\end{proof}

Moreover, we can show that actually the fixed points of $\mathcal{F}$ are weak solutions. The next proposition proves the equivalence between the weak formulation and our mild formulation.

   \begin{lemma}[Mild implies weak]\label{lem:MildImpliesWeak}
Let $\mu\in X_T$ be a fixed point of $\mathcal{F}$. Then for every $f\in\mathcal E$ the map $t\mapsto\langle\mu_t,f\rangle$
is absolutely continuous and satisfies the weak formulation of the coagulation transport equation in the sense of \eqref{eq:weakSolution}.
\end{lemma}
\begin{proof}
Fix $f\in\mathcal E$. Then,
\begin{equation}
\label{eq:mildFormulationTestPairing} \pint{\mu_t,f} =\pint{S_{\Phi}(t)\nu,f}+\int_0^t\pint{S_\Phi(t-s)Q(s,\mu_s),f} \,ds.
\end{equation}

For each $t>0$, the map $s\mapsto \pint{S_\Phi(t-s)Q(s,\mu_s),f}$ is bounded and continuous and also, for fixed $s<t$, by Proposition \ref{prop:weakContinuityTransport} the map $t\mapsto \pint{S_\Phi(t-s)Q(s,\mu_s),f}$ is differentiable with 
\[
\frac{d}{dt}\pint{S_\Phi(t-s)Q(s,\mu_s),f}=-\pint{S_\Phi(t-s)Q(s,\mu_s),\psi f'}. 
\]
Moreover, $s\mapsto\frac{d}{dt}\pint{S_\Phi(t-s)Q(s,\mu_s),f}$ is continuous and bounded on $[0,T]$ since 
$\abs{\pint{S_\Phi(t-s)Q(s,\mu_s),\psi f'}}\leq \abs{\psi f'}_\infty\frac{\norm{\mu_t}_\infty(\norm{\mu_t}_\infty+1)}{s+\kappa}$. Then, by dominated convergence theorem we can differentiate under the integral sign and apply Leibniz rule obtaining 
\[
\frac{d}{dt} \int_0^t\pint{S_\Phi(t-s)Q(s,\mu_s),f} \,ds = \pint{Q(t,\mu_t),f}-\int_0^t
\pint{S_\Phi(t-s)Q(s,\mu_s),\psi f'}
\,ds. 
\]
Then, differentiating \eqref{eq:mildFormulationTestPairing} and rearranging the terms we get
\begin{align*}
   \frac{d}{dt}\pint{\mu_t,f}&=-\pint{S_{\Phi}(t)\nu,\psi f'} +\pint{Q(t,\mu_t),f}-\int_0^t
\pint{S_\Phi(t-s)Q(s,\mu_s),\psi f'}
\,ds  \\
&= \pint{Q(t,\mu_t),f}-\pint{\mu_t,\psi f'}.
\end{align*}
This completes the proof.
\end{proof}

We are ready to use Banach fixed point theorem which will act as the main tool to show local existence and uniqueness of solutions. To use the aforementioned theorem in this context we need to show two things: that $\mathcal{F}$ maps a closed ball of $X_T$ into itself and that $\mathcal{F}$ is a contraction restricted to that closed subspace. So, let $R>0$ and let us denote the closed ball of radius $R$ in $X_T$ as $B_R$, that is, $B_R=\set{\mu_t\in X_T\,: \,\norm{\mu_t}_\infty\leq R}.$ We will start by showing that in this closed subspace of $X_T$ the operator $F$ is a contraction.

\begin{proposition}  If $T$ is sufficiently small, then $\mathcal{F}(B_R)\subset B_R$. Also, for $\mu_t, \tilde{\mu}_t\in B_R$, we have the bound   
    \[
    \norm{\mathcal{F}(\mu_t)-\mathcal{F}(\tilde{\mu}_t)}_{\infty}\leq \left[ (2R+1)\log\,\frac{T+\kappa}{\kappa}\right]\norm{\mu_t-\tilde{\mu}_t}_{\infty}.
    \]
   which implies that a sufficient condition for $\mathcal{F}$ to be a contraction is that 
\begin{equation}\label{eq:secondSmallnessCondition}
        (2R+1)\log \frac{T+\kappa}{\kappa}<1.
    \end{equation}
\end{proposition}
\begin{proof} 
Applying \eqref{eq:pettisIntegralBound} we obtain,
\begin{align*}
    \norm{\mathcal{F}(\mu_t)}_{\infty}&\leq \norm{\nu}+\norm{\mathcal{I}(\mu_t)}_\infty\\
    &\leq \norm{\nu}
    +(\norm{\mu_t}_\infty+\norm{\mu_t}_\infty^2)\log\left(\frac{T+\kappa}{\kappa}\right)\\
    &\leq \norm{\nu}+R(R+1)\log\,\frac{T+\kappa}{\kappa}.
\end{align*}
Thus, a sufficient condition for $F$ to map $B_R$ into itself is 
\begin{equation}\label{eq:firstSmallnessConditon}
\norm{\nu}+R(R+1)\log\,\frac{T+\kappa}{\kappa}\leq R.
\end{equation}

Now, let $\mu_t, \tilde{\mu}_t\in B_R$, using \eqref{eq:CoagulationLipschitzBound} we obtain
    \[
    \norm{Q(s, \mu_t) - Q(s, \tilde{\mu}_t)}\leq \frac{2R+1}{s+\kappa}\norm{\mu_t - \tilde{\mu}_t}_{\infty}.
    \]
    Thus,
    \[
    \norm{\mathcal{F}(\mu_t) - \mathcal{F}( \tilde{\mu}_t)}\leq (2R+1)\norm{\mu_t -  \tilde{\mu}_t}_{\infty}\int_0^t\frac{1}{s+\kappa}\,ds,
    \]
    which gives us the desired bound. 

\end{proof}
  
The following lemma will show that the operator $\mathcal{F}$ can be defined in a way in which it is evident that the positivity of the measures is preserved, even though we do not change the definition per se, but we will use the property to prove positivity preservation using a contradiction argument.
 \begin{lemma}\label{lem:additiveForm}
   Let $\mu_t\in X_T$ a fixed point of $\mathcal{F}$. Then, for every $f\in C_b(\R^+)$
   \begin{equation}
\label{eq:aditiveMildFormula}\pint{\mu_t,f}=\frac{\kappa}{t+\kappa}\pint{S_\Phi(t)\nu,f}+\frac{1}{t+\kappa}\int_0^t \pint{S_\Phi(t-s)\mu_s\star\mu_s,f}\, ds
   \end{equation}
\end{lemma}
\begin{proof}
Let $f\in C_b(\R^+)$. To analyze how the map $t\mapsto\pint{\mu_t,f}$ evolves up to a certain time, let us fix $\tilde t\in[0,T]$. We consider the pairing of $\mu_t$ with the auxiliary function $f\circ\Phi_{\tilde t-t}$. Thus, using the semigroup property \eqref{eq:semigroupProperties} we see that  
\begin{align*}
    \pint{\mu_t,f\circ\Phi_{\tilde t-t}} &=\pint{S_{\Phi}(t)\nu,f\circ\Phi_{\tilde t-t}}+\int_0^t\pint{S_\Phi(t-s)Q(s,\mu_s),f\circ\Phi_{\tilde t-t}} \,ds\\
&=\pint{\nu,f\circ\Phi_{\tilde t} }+\int_0^t \frac{1}{s+\kappa}\left[\pint{\mu_s \star\mu_s,f\circ\Phi_{\tilde{t}-s} }-\pint{\mu_s ,f\circ\Phi_{\tilde{t}-s} }\right]\,ds
\end{align*}
Note that in this last equation $\pint{\nu,f\circ\Phi_{\tilde t} } $ is a constant. Let us define then $G(t):=\pint{\mu_t,f\circ\Phi_{\tilde t-t}}$ and $R(s):=\pint{\mu_s \star\mu_s,f\circ\Phi_{\tilde{t}-s} }$. Thus, we obtain the following Volterra integral equation for $G$:
\[
G(t)=\pint{\nu,f\circ\Phi_{\tilde t} }+\int_0^t \frac{1}{s+\kappa}R(s)\,ds - \int_0^t \frac{1}{s+\kappa}G(s)\,ds
\]
Since $R$ and $G$ are continuous, we conclude that $G$ is absolutely continuous and therefore differentiable almost everywhere in $t\in(0,\tilde t\, )$ obtaining the equation
\[
\frac{d}{dt}G(t)+\frac{1}{t+\kappa}G(t)=\frac{1}{t+\kappa}R(t),
\]
which we can solve in the following way by using an integrating factor
\[
\frac{d}{dt}\left[\frac{\kappa+t}{\kappa}G(t)\right]=\frac{1}{\kappa} R(t).
\]
Integrating from $0$ to $\tilde t$ and rearranging the terms we get
\[
G(\tilde t\,)=\frac{\kappa}{\tilde t+\kappa}G(0)+\frac{1}{\tilde t+\kappa}\int_0^{\tilde t}R(s)\, ds,
\] which, in the context of the measure, is 
\[
\pint{\mu_{\tilde t},f}= \frac{\kappa}{\tilde t+\kappa}\pint{\nu,f\circ\Phi_{\tilde t}}+\frac{1}{\tilde t+\kappa}\int_0^{\tilde t}\pint{\mu_s \star\mu_s,f\circ\Phi_{\tilde{t}-s} }\, ds.
\]
This finishes the proof.
\end{proof}

Now, with this additive description of the fixed point we can prove that the operator preserves positivity, that is, if the negative part of the  initial condition is null, then all the evolution of the measure solution will remain non-negative, as a measure. We give an argument using duality.

\begin{proposition}[Positivity]
\label{prop:positivityOfFixedPoints}
Let $\nu\in \mathcal{M}(\R^+)$ be a positive measure and $\mu_t\in X_T$ a fixed point of $\mathcal{F}$. Then, $\mu_t$ is a positive measure for all $t\in[0,T]$ 
\end{proposition}
\begin{proof}
We will argue by contradiction assuming that the set $\mathcal{A}=\set{t\in[0,T]\, :\,  \mu_t^->0}$ is non empty. Let $\tilde t=\inf\mathcal{A}$. Since $\mu_0=\nu\geq0$, we know that $\mu_s$ is a positive measure for all $s<\tilde t$. Now, let us take $f\in C_c(\R^+)$ such that $f\geq 0$. Then, applying \eqref{eq:aditiveMildFormula} we see that  $\pint{\mu_{\tilde t},f}\geq 0$, since all the integrals involve positive measures and positive functions. This forces $\mu_{\tilde t }$ to be a positive measure, implying that $\tilde t \not\in A$.

Let $\mu_t=\mu_t^+-\mu_t^-$ the Jordan decomposition of $\mu_t$ for $t>\tilde t$. For $s\in[0,T]$ we have 
\[\mu_s \star \mu_s = \left(\mu_s^+ \star \mu_s^+ + \mu_s^- \star \mu_s^-\right) - 2\left(\mu_s^+ \star \mu_s^-\right).\]
Since the transport operator preserves positivity and $\nu$ is a positive measure, formula \eqref{eq:aditiveMildFormula} gives us a decomposition of $\mu_t$ as the difference of non-negative measures as 
$\mu_t=\eta_t^+-\eta_t^-$ where the measures $\eta_t^+$ and $\eta_t^-$ act on $f\in C_b(\R)$ as follows:
\[
\pint{\eta_t^+,f}= \frac{\kappa}{ t+\kappa}\pint{\nu,f\circ\Phi_{t}}+\frac{1}{ t+\kappa}\int_0^{ t}\pint{\mu_s^+ \star \mu_s^+ + \mu_s^- \star \mu_s^-,f\circ\Phi_{{t}-s} }\, ds.
\]
and
\[
\pint{\eta_t^-,f}= \frac{2}{t+\kappa}\int_0^{ t}\pint{\mu_s^+ \star \mu_s^-,f\circ\Phi_{{t}-s} }\, ds.
\]
The minimality of the Jordan decomposition implies that $\mu_t^-\leq \eta_t^-$ and, this integrating against the function $f\equiv1$ we get the estimation 
\[
\norm{\mu_t^-}\leq\frac{2}{t+\kappa}\int_0^t \norm{\mu_s^-}\norm{\mu_s^+}\ ds
\]
Using that $\mu_s^-=0$ for $t\leq\tilde t$ 
and the bound $\norm{\mu_s^+}\leq \norm{\mu_t}_\infty$ we can simplify the bound as 
\[
\norm{\mu_t^-}\leq\frac{2 \norm{\mu_t}_\infty}{\kappa}\int_{\tilde t}^t \norm{\mu_s^-}\ ds
\]
Let $\varepsilon>0$. By taking the supremum of $\norm{\mu_t^-}$ for $t\in[\tilde t,\min\set{T,\tilde t+\varepsilon}]$ we obtain 
\[
\sup_{t\in[\tilde t,\min\set{T,\tilde t+\varepsilon}]}\norm{\mu_t^-}\leq \frac{2 \norm{\mu_t}_\infty}{\kappa}\varepsilon \sup_{t\in[\tilde t,\min\set{T,\tilde t+\varepsilon}]}\norm{\mu_t^-}.
\]
If $\varepsilon$ is sufficiently small, then $\norm{\mu_t^-}=0$ for $t\in(\tilde t, \tilde t+\varepsilon)$, contradicting the minimality of $\tilde t$.
\end{proof}
Finally, we will apply the following standard result to show that the local solution can be extended to an arbitrary time. Though it is a standard result in PDE analysis result, we include it for the sake of completeness. 
  \begin{lemma}  \label{lem:timeExtention}
       Let $T_{\max}:=\sup\{T>0:\text{ there exists a unique solution }\mu_t\in X_T\}$ be the maximal existence time of the solution constructed by the fixed point method. Then exactly one of the following holds:
\begin{enumerate}
  \item $T_{\max}=+\infty$ (global existence);
  \item $T_{\max}<\infty$ and $\displaystyle\lim_{t\to T_{\max}^-}
  \norm{\mu_t}
=+\infty$ (finite time blow-up).
\end{enumerate}
Equivalently, if $\sup_{t\in[0,T_{\max})}\norm{\mu_t}<\infty$ then $T_{\max}=+\infty$.
  \end{lemma}
  \begin{proof}
      Suppose, toward a contradiction, that $T_{\max}<\infty$ but the norm remains bounded on $[0,T_{\max})$. Then there exists $R_*<\infty$ such that
\[
\norm{\mu_t} \le R_* \qquad\text{for all }t\in[0,T_{\max}).
\]

There exists a time $T_{R_*}>0$ (depending only on $R_*$) such that the fixed point map $F$ is a contraction on the ball $B_{R_*}$ when considered on any interval of length at most $T_{R_*}$. The construction of the fixed point is local in time and time-translation invariant in the following sense: if $\tilde\mu$ is the solution on $(0,t_0)$ then, treating $\tilde\mu_{t_0}$ as new initial data, the same fixed point argument produces a unique solution on $(t_0,t_0+T_{R_*})$ provided 
$\norm{\tilde\mu_{t_0}}\leq R_*$.

Since $\norm{\mu_{t_0}}\le R_*$ for every $t_0<T_{\max}$, we may pick $t_0\in(0,T_{\max})$ so close to $T_{\max}$ that
\[
t_0 + T_{R_*} > T_{\max}.
\]
Applying the fixed point argument with initial time $t_0$ and initial data $\mu_{t_0}$ yields a unique solution on $(t_0,t_0+T_{R_*})$. By uniqueness this solution coincides with the original solution on $(0,t_0)$, and hence provides an extension of the original solution beyond $T_{\max}$, contradicting the definition of $T_{\max}$ as the supremum of existence times.

Therefore the assumption that $T_{\max}<\infty$ while $\sup_{t<T_{\max}}\norm{\mu_t}<\infty$ is impossible. This proves that either $T_{\max}=\infty$ or $\limsup_{t\to T_{\max}^-}\norm{\mu_t}=+\infty$, which is the stated blow-up alternative.
\end{proof}

We can now write down the proof of the main theorem which will be obtained as a result of all the theory developed throughout this section. \begin{proof}[Proof of Theorem \ref{thm:mainTheorem}] Let us fix an initial condition $\nu\in\mathcal{M}_P(\R^+)$.
 First we note that the constants $R>0$ and $T>0$ can be chosen so that $\mathcal{F}$ maps $B_R$ into itself and is a contraction in that subspace: We only need to verify that it is possible to choose $R>0$ and $T>0$ such that both \eqref{eq:firstSmallnessConditon} and \eqref{eq:secondSmallnessCondition} hold. This system can be rewritten as 
\[
\log \frac{T+\kappa}{\kappa}<\frac{1}{2R+1}\quad\text{and}\quad\log \frac{T+\kappa}{\kappa}\leq\frac{R-\norm{
\nu}}{R(R+1)}.
\]
Then it is clear that we can choose $R>\norm{
\nu}$ and $T$ small enough so that both inequalities hold. Under such conditions $\mathcal{F}$ has a unique fixed point in $B_R$.

Let us choose a $T,R>0$ such that there is a unique fixed point $\mu_t\in X_T$ of $\mathcal{F}$ with initial condition $\nu$. Then Proposition \ref{prop:positivityOfFixedPoints} implies that this solution remains positive until time $T$, Lemma \ref{lem:MildImpliesWeak} implies that this fixed point is a local weak solution, Lemma \ref{lem:probabilityConcervation}  implies that the mass is conserved, for this reason we can apply Lemma \ref{lem:timeExtention} to obtain that the unique local weak solution can be extended up to an arbitrary time.
\end{proof}

\section{Regularity and formula for a solution starting with a Dirac measure}
\label{sec:Regularity}
In this section we deal with some regularity conditions and how dynamics naturally propagates them from the initial condition. We will deal with the moment preservation, absolute continuity with respect to Lebesgue measure and a formula for the solution in a case relevant from the modeling standpoint. The moment preservation can be easily derived as an \textit{a priori} estimate:

\begin{proposition} \label{prop:momentBound}
    Let $\mu_t:[0,T]\to \mathcal{M}(\R^+)$ be a local weak solution with initial condition $\nu$. If $\mu_t$ is a positive measure for $t\in[0,T]$ and $\nu$ is a probability measure, then the following estimate holds
 \[\int_{\mathbb R^+} x\,\mu_t(dx)\leq \frac{t+\kappa}{\kappa}\intRp x\,\nu(dx),\quad\text{for all }t\leq T.
\]
In particular, if $\nu$ has a finite first moment, then $\mu_t$ has a finite first moment for every $ 0\leq t\leq T$.
    
\end{proposition}
\begin{proof}
    To approximate the identity on $\R^+$ we will construct a suitable family of functions in our test space $\mathcal{E}$ allowing us to pass to the limit using the monotone convergence theorem. to this end we choose non-decreasing and concave cut function  $\varphi \in C^1(\R^+)$ such that $\varphi(x)=x$ if $x\in[0,1]$ and $\varphi$ is constant in $[2,\infty)$. On the transition interval $(1,2)$, it smoothly interpolates between these regimes. Then, for $R>0$ we will consider the following family of functions $f_R(x):=R\varphi(x/R)$. By construction $f_R$ is bounded and $f_R'$ has support in $[0,2R]$, thus $\psi f_R'$ is bounded. So, $f_R\in\mathcal{E}$ and converges monotonically to the identity function as $R\to \infty$. 
    Set
\[M_R(t):=\langle\mu_t,f_R\rangle.\] Using $f_R$ in \eqref{eq:weakSolution} we obtain
\[
M_R(t)=M_R(0)-\int_0^t
\pint{\mu_s,\psi f'_R}\,ds+\int_0^t\frac{1}{s+\kappa}\Big(\pint{\mu_s\star\mu_s,f_R}-\pint{\mu_s,f_R}\Big)\,ds.
\]
Even more, we obtain 
\begin{equation}\label{eq:momentEstimationR}
     \frac{d}{dt}M_R(t)=-\pint{\mu_t,\psi f'_R}+\frac{1}{t+\kappa}\Big(\pint{\mu_t\star\mu_t,f_R}-M_R(t)\Big).
\end{equation}
Since $f_R$ is concave and $f_R(0)\geq0$, we know that $f_R$ is sub-additive. Thus,
\begin{align*}
\langle\mu_s\star\mu_s,f_R\rangle&=\intRp \intRp  f_R(x+y)\,\mu_s(dx)\,\mu_s(dy)\\
    &\le\intRp \intRp(f_R(x)+f_R(y))\,\mu_s(dx)\,\mu_s(dy)=2M_R(s).
\end{align*} Therefore \eqref{eq:momentEstimationR} becomes
\begin{equation}\label{eq:completeMomentBound}
      \frac{d}{dt}M_R(t)\leq -\pint{\mu_t,\psi f'_R}+\frac{1}{t+\kappa}M_R(t).
\end{equation}
Inasmuch as $\langle\mu_s,\psi f'_R\rangle\geq0 $ for all $s\geq0$, we can simplify even more 
\[
\frac{d}{dt}M_R(t)\leq \frac{1}{t+\kappa}M_R(t).
\]
This bound is enough to solve the Grönwall type inequality and obtain
\[
M_R(t)\leq M_R(0)\frac{t+\kappa}{\kappa}.
\] Since $f_R\uparrow x$ point-wise and $\mu_t$ is finite, monotone convergence gives $M_R(t)\xrightarrow[R\to\infty]{}\intRp x\,\mu_t(dx)$. Because $\int_{\mathbb R^+} x\,\nu(dx)<\infty$, we pass to the limit to obtain the desired bound.
\end{proof}

In the previous proof one can choose for instance the mollifier 
\[
\varphi(x)=\begin{cases}
    x\quad &\text{if } x\in[0,1],\\
    x-(x-1)^3+\frac{1}{2}(x-1)^4&\text{if } x\in[1,2],\\
    \frac{3}{2}&\text{if } x\in[2,\infty).
\end{cases}
\]
 This function satisfies all the requirements needed in the proof.

To explore the remaining properties we will return to the mild formulation, more precisely, to the modification of the fixed point operator we used to prove the invariance of the positive cone under the dynamics induced by the transport and coagulation in Proposition \ref{prop:positivityOfFixedPoints}. Let us consider the following operator on $X_T$ defined weakly as 
\begin{equation}
\label{eq:additiveOperatorForm}
    \pint{\mathcal{G}(\mu_t),f}:= \frac{\kappa}{t+\kappa}\pint{S_\Phi(t)\nu,f}+\frac{1}{t+\kappa}\int_0^t \pint{S_\Phi(t-s)\mu_s\star\mu_s,f}\, ds \quad\text{for all}\quad f\in C_0(\R^+).
\end{equation}
This is obtained after applying an integration factor to $\mathcal{F}$. Moreover, Lemma \ref{lem:additiveForm} shows that fixed points of 
$\mathcal{F}$ are fixed points of $\mathcal{G}$, and performing exactly the same computations in reverse order we can see that in fact the fixed points are the same. We include the argument for the sake of completeness.
\begin{proposition}
Let $\mu_t\in X_T$. If $\mathcal{G}(\mu_t)=\mu_t$ for all $t\in[0,T]$, then $\mathcal{F}(\mu_t)=\mu_t$ for all $t\in[0,T]$.
\end{proposition}
\begin{proof}
   Let $f\in C_0(\R^+)$ and fix $t\in[0,T]$. Then \eqref{eq:additiveOperatorForm} is equivalent to 
   \[
   \int_0^t \frac{d}{ds}\left[\frac{s+\kappa}{\kappa}\pint{\mu_s,f\circ\Phi_{t-s}}\right]=\frac{1}{\kappa}\int_0^t \pint{S_\Phi(t-s)\mu_s\star\mu_s,f}\, ds.
   \]
   Then, 
   \[
   \frac{d}{ds}\left[\frac{s+\kappa}{\kappa}\pint{\mu_s,f\circ\Phi_{t-s}}\right]=\frac{1}{\kappa} \pint{S_\Phi(t-s)\mu_s\star\mu_s,f}.
   \] 
   Expanding the derivative using the product rule and rearranging the terms we obtain
   \[
   \frac{d}{ds}\pint{\mu_s,f\circ\Phi_{t-s}}=\frac{1}{s+\kappa}\left(\pint{S_\Phi(t-s)\mu_s\star\mu_s,f}-\pint{S_\Phi(t-s)\mu_s,f}\right),
   \]
   which we can integrate from $0$ to $t$ to conclude that $\mathcal{F}(\mu_t)=\mu_t$, since $f\in C_0(\R^+)$ was arbitrary.
\end{proof}

 Similar to the proof of Theorem \ref{thm:mainTheorem} at the end of Section \ref{sec:MildFormulation}, we have that the operator  $\mathcal{G}$ maps a closed ball $B_R$ of $X_T$ to itself and is a contraction if the following both conditions hold:  $T\leq\frac{\kappa(R-\norm{\nu})}{R^2}$ and $T<\frac{\kappa}{R^2}$. As a consequence, $\mathcal{G}$ has a unique fixed point.
We will take advantage of the convergence to the fixed point in total variation norm using an iterative argument.

\begin{proposition} \label{prop:lebesgueRegularity}
   Let $\mu_t$ be a weak solution on a finite interval $[0,T]$. If $\nu\ll\mathcal{L}$, then $\mu_t\ll\mathcal{L}$ for all $t\in[0,T]$.
\end{proposition}
\begin{proof}
    We know that
     $\mu_t$ is a fixed point of $\mathcal{G}$ in $X_T$ for a possibly small $T$. So, let us consider the following iterative scheme, let us define recursively on $n\in\mathbb{N}$ the measures: $\mu_t^{(0)}=\frac{\kappa}{t+\kappa}S_{\Phi}(t)\nu$ and $\mu_t^{(n+1)}:=\mathcal{G}(\mu_t^{(n)})$.

     We will proceed by mathematical induction on the iterative step to show that every iterated measure is absolutely continuous. Since $\nu\ll \mathcal{L}$, we can say that $\nu=u_0\, dx$ for some non-negative function $u_0\in L^1(\R^+)$. Let $f\in C_0(\R^+)$, then by the change of variable $y=\Phi_t(x)$, we get
     \[
\pint{\mu_t^{0},f}=\frac{\kappa}{t+\kappa}\int_{\Phi_t(\R^+)} f(y)u_0(\Phi_{-t}(y))\frac{d}{dy}\Phi_{-t}(y)\, dy.
     \]
     Since $\frac{d}{dy}\Phi_{-t}(y)=\op{exp}\left(\int_0^t\psi'(\Phi_{-s}(y))\right)\,ds$, and the weak derivative of $\psi$ is essentially bounded, then the change of variable is justified and the function 
     \[
     u^{(0)}(t,y):=\frac{\kappa}{t+\kappa}\I{\Phi_t(\R^+)}(y)u_0(\Phi_{-t}(y))\op{exp}\left(\int_0^t\psi'(\Phi_{-r}(y))\,dr\right)
     \] is a legitimate $L^1(\R^+)$ density for each fixed $t$.

     Now, lets assume inductively that $\mu_t^{(n)}$ has a density $u^{(n)}(t,y)$. Since the convolution of two absolutely continuous measures (with respect to the Lebesgue measure) is absolute continuous and its density is exactly the convolution of their densities, we know that $\mu_s\star\mu_s$ has density given by 
     \[
u_C(s,y):=\int_0^yu^{(n)}(s,y-z)u^{(n)}(s,z)\, dz.
     \]
     Then, using the same technique as above and interchanging the order of the integrals by Fubini's theorem we see that $\mu_t^{(n+1)}$ has density 
     \[
     u^{(n+1)}(t,y):=u^{(0)}(t,y)+\frac{1}{t+\kappa}\int_0^t\I{\Phi_{t-s}(\R^+)}(y)u_C(s,\Phi_{-(t-s)}(y))\op{exp}\left(\int_0^{t-s}\psi'(\Phi_{-r}(y))\,dr\right)\,ds.
     \]
This proves that for all $n\in\mathbb{N}$, $\mu_t^{(n)}\ll\mathcal{L}$. 
     Finally, since $\mathcal{G}$ is a contraction for $T$ small enough, then, for every $t\leq T$ $\mu_t^{(n)}$ converges in total variation to the solution $\mu_t$. Since the subspace of absolutely continuous measure is closed in this topology, we conclude that indeed $\mu_t\ll\mathcal{L}$. The mild formulation \eqref{eq:additiveOperatorForm} guarantees that we can extend for arbitrary $T$.
\end{proof}

Given $x_0>0 $, we can use same technique as above to study the problem 
\begin{equation}\label{eq:InitialDelta}
    \begin{cases}\displaystyle
\partial_t\mu_t=\partial_x(\frac{x^2}{2} \mu_t)+\frac{1}{t+\kappa}(\mu_t\star\mu_t-\mu_t),\\
        \mu_0=\boldsymbol{\delta}_{x_0},
    \end{cases}
\end{equation}
This problem is relevant from the modeling point of view since it describes the evolution of a monodisperse population where all the species have a rescaled initial mass $x_0$. Again, we will use the iteration scheme to prove the following result. 
\begin{proposition}\label{prop:DeltaSolutionFormula}
Let $\mu_t$ be the weak solution of \eqref{eq:InitialDelta} in $[0,T]$.
There exists a function $u:[0,T]\times \R^+\to\R^+$ such that  $u(t,\cdot)\in L^1(\R^+)$ for all $t\in[0,T]$ and 
\begin{equation}
\label{eq:DeltaSolutionLebesgueDecomposition}\mu_t(dx)=\frac{\kappa}{t+\kappa}\boldsymbol{\delta}_{\Phi_t(x_0)}+u(t,x)dx.
\end{equation}
\end{proposition}
\begin{proof}
    In this case the flow has an explicit formula representation, that is $\Phi_t(x)=\frac{2x}{2+tx}.$ Let us consider again the iterative scheme: $\mu_{t}^{(0)}:=\frac{\kappa}{t+\kappa}\boldsymbol{\delta}_{\frac{2x_0}{2+tx_0}}$ and  $\mu_{t}^{(n+1)}:=\mathcal{G}(\mu_{t}^{n})$ for $n\in\mathbb{N}$ greater than $0$. 

    We claim that for each $n\geq 1$ the measure $\mu_t^{(n)}$ has a single atomic part given by $\mu_t^{(0)}$ and the rest is absolutely continuous. First let us do the computation for the case $n=1$: Let $f\in C_0(\R^+)$, then, since the convolution of two Dirac measures is another Dirac measure with the sum of the center points, we will have 
    \[
\pint{\mu_t^{(1)},f}=\frac{\kappa}{t+\kappa}f(\Phi_t(x_0))+\frac{1}{t+\kappa}\int_0^t\frac{\kappa^2}{(s+\kappa)^2}f(\Phi_{t-s}(2\Phi_s(x_0)))\, ds.
    \]
Note that in the integrand $f$ varies smoothly as $s$ varies between $0$ and $t$, this defines a density as we shall see. After a change of variable in in the integral we obtain 
\[
\pint{\mu_t^{(1)},f}=\frac{\kappa}{t+\kappa}f(\Phi_t(x_0))+\frac{1}{t+\kappa}\int_{\Phi_t(2x_0)}^{2\Phi_t(x_0)}\frac{4\kappa^2}{((\kappa+2/x_0+2t)y-4)^2}f(y)\, dy.
\]
Since $f$ was arbitrary, this proves that $\mu_t^{(1)}$ has a decomposition $\mu_t^{(1)}=\frac{\kappa}{t+\kappa}\boldsymbol{\delta}_{\Phi_t(x_0)}+u^{(1)}(t,x)dx$, where 
\[u^{(1)}(t,x)=\frac{1}{t+\kappa}\I{ [\Phi_t(2x_0),2\Phi_t(x_0)]}(x)\frac{4\kappa^2}{((\kappa+2/x_0+2t)x-4)^2}.\]
    Now, lets inductively assume that $\mu_t^{(n)}$ has a decomposition $\mu_t^{(n)}=\frac{\kappa}{t+\kappa}\boldsymbol{\delta}_{\Phi_t(x_0)}+u^{(n)}(t,x)dx$ for a given $u^{(n)}(t,\cdot)\in L^1(\R^+)$ and let $f\in C_0(\R^+)$. The convolution $\mu_{s}^{(n)}\star\mu_{s}^{(n)}$ yields 
    \[
\mu_{s}^{(n)}\star\mu_{s}^{(n)}=\frac{\kappa^2}{(s+\kappa)^2}\boldsymbol{\delta}_{2\Phi_s(x_0)}+2\frac{\kappa}{s+\kappa}\boldsymbol{\delta}_{\Phi_s(x_0)}\star u^{(n)}(s,dx)+ (u^{(n)}(s,dx))\star (u^{(n)}(s,dx))
    \]
As we saw in the case $\mu_t^{(1)}$, the first term on the right will lead to a density. We know that the measure $(u^{(n)}(s,dx))\star (u^{(n)}(s,dx))$ is absolutely continuous and we will call its density $u^{(n)}_C$ just for the sake of simplicity. On the other hand, the measure $2\frac{\kappa}{s+\kappa}\boldsymbol{\delta}_{\Phi_s(x_0)}\star u^{(n)}(s,dx)$ will also have density, since the convolution with the Dirac measure acts as a translation; its density is $
u^{(n)}_T(s,y):=2\frac{\kappa}{s+\kappa}\I{[\Phi_s(x_0),\infty)}(y)u^{(n)}(s,y-\Phi_s(x_0))$. Thus,
\begin{align*}
    \pint{\mu_t^{(n+1)},f}=&\frac{\kappa}{t+\kappa}f(\Phi_t(x_0))+\intRp f(x)u^{(1)}(t,x)\,dx\\&+\frac{1}{t+\kappa}\int_0^t \intRp f(\Phi_{t-s}(y))(u^{(n)}_T(s,y)+u^{(n)}_C(s,y))\, dyds
\end{align*}
Changing variables and applying Fubini's theorem we conclude that $\mu_t^{(n+1)}=\frac{\kappa}{t+\kappa}\boldsymbol{\delta}_{\Phi_t(x_0)}+u^{(n+1)}(t,y)dy$ with
\[
u^{(n+1)}(t,y):=u^{(1)}(t,y)+\frac{1}{t+\kappa}\int_0^t\I{\Phi_{t-s}(\R^+)}(y)u_{TC}(s,\Phi_{-(t-s)}(y))\op{exp}\left(\int_0^{t-s}\psi'(\Phi_{-r}(y))\,dr\right)\,ds,
\]
 where  $u_{TC}=u_{T}+u_{C}$.

Finally, since the convergence in total variation preserves the Lebesgue decomposition we conclude that, for a possible small $T$, $\mu_t$ has the desired form \eqref{eq:DeltaSolutionLebesgueDecomposition} for all $t<T$. Again, the mild formulation \eqref{eq:additiveOperatorForm} guaranties that we can extend for arbitrary $T$.
\end{proof}

The above results constitute the content of Theorem \ref{thm:secondMainThm} which we can  prove now. 

\begin{proof}[Proof of Theorem \ref{thm:secondMainThm}.] The assertion follows immediately by assembling the results developed in this section: Part $(i)$ is stated in Proposition \ref{prop:momentBound}, $(ii)$ is stated in Proposition \ref{prop:lebesgueRegularity} and, finally $(iii)$ is stated in Proposition \ref{prop:DeltaSolutionFormula}.
\end{proof}

\section*{Acknowledgments} EAA acknowledges partial support from SECIHTI grant for doctoral studies. ASJ acknowledges DGAPA-PAPIIT-UNAM grant IN-105726 for partial support. LFLR was partially supported by DGAPA-UNAM, program PAPIIT, grant IN113225.

\printbibliography

\end{document}